# Accurate and efficient explicit approximations of the Colebrook flow friction equation based on the Wright-Omega function

Dejan Brkić [1,2,*] and Pavel Praks [1,3,*]

[1] European Commission, Joint Research Centre (JRC), Directorate C: Energy, Transport and Climate, Unit C3: Energy Security, Distribution and Markets, Via Enrico Fermi 2749, 21027 Ispra (VA), Italy

[2] Alfatec, Bulevar Nikole Tesle 63/5, 18000 Niš, Serbia

[3] IT4Innovations National Supercomputing Center, VŠB—Technical University of Ostrava, 17. listopadu 2172/15, 708 00 Ostrava, Czech Republic

[*] Correspondence:

dejanbrkic0611@gmail.com (D.B.); ORCID id: 0000-0002-2502-0601

pavel.praks@ec.europa.eu or pavel.praks@gmail.com (P.P.); ORCID id: 0000-0002-3913-7800

**Abstract:** The Colebrook equation is a popular model for estimating friction loss coefficients in water and gas pipes. The model is implicit in the unknown flow friction factor $f$. To date, the captured flow friction factor $f$ can be extracted from the logarithmic form analytically only in the term of the Lambert $W$-function. The purpose of this study is to find an accurate and computationally efficient solution based on the shifted Lambert $W$-function also known as the Wright $\Omega$-function. The Wright $\Omega$-function is more suitable because it overcomes the problem with the overflow error by switching the fast growing term $y = W(e^x)$ of the Lambert $W$-function to the series expansions that further can be easily evaluated in computers without causing overflow run-time errors. Although the Colebrook equation transformed through the Lambert $W$-function is identical to the original expression in term of accuracy, a further evaluation of the Lambert $W$-function can be only approximate. Very accurate explicit approximations of the Colebrook equation that contains only one or two logarithms are shown. The final result is an accurate explicit approximation of the Colebrook equation with the relative error of no more than 0.0096%. The presented approximations are in the form suitable for everyday engineering use, they are both accurate and computationally efficient.

## 1. Introduction

The Colebrook equation; Eq. (1), is an empirical relation which in its native form relates implicitly the unknown Darcy's flow friction factor $f$ with the known Reynolds number $R$ and the known relative roughness of inner pipe surface $\varepsilon^*$ [1,2]. Engineers use it at defined domains of the input parameters: $4000 < R < 10^8$ and for $0 < \varepsilon^* < 0.05$. The Colebrook equation is transcendental (cannot be expressed in the term of elementary functions), the implicitly given function in respect to the unknown flow friction factor $f$

$$\frac{1}{\sqrt{f}} = -2 \cdot \log_{10}\left(\frac{2.51}{R} \cdot \frac{1}{\sqrt{f}} + \frac{\varepsilon^*}{3.71}\right), \tag{1}$$

The Colebrook equation; Eq. (1) has also an exact explicit analytical form in the term of the Lambert $W$-function; Eq. (2) [3,4] that is also transcendental, but which can be evaluated through the numerous thoroughly tested procedures of various accuracy and complexity developed for various applications in physics and engineering [5].

$$\left.\begin{array}{l}\frac{1}{\sqrt{f}} = \frac{2}{\ln(10)} \cdot \left(\ln\left(\frac{R}{2.51} \cdot \frac{\ln(10)}{2}\right) + W(e^x) - x\right) \\ x = \ln\left(\frac{R}{2.51} \cdot \frac{\ln(10)}{2}\right) + \frac{R \cdot \varepsilon^*}{2.51 \cdot 3.71} \cdot \frac{\ln(10)}{2}\end{array}\right\}, \tag{2}$$

The parameter $x$ in Eq. (2) depends on the input parameters; the Reynolds number $R$ and the relative roughness of inner pipe surface $\varepsilon^*$. Its domain is 7.51<$x$<618187.84. The Lambert $W$-based Colebrook equation; Eq. (2), contains the fast growing term $W(e^x)$, which cannot be stored in computer registers due to the runtime overflow error for the certain combinations of the Reynolds number $R$ and the relative roughness of inner pipe surface $\varepsilon^*$ that can easily occur in everyday engineering practice [6,7]. The problem can be solved using the Wright $\Omega$-function, a cognate of the Lambert $W$-function, which uses a shifted, not fast-growing argument [4,8,9].

This paper presents few approximate solutions of the transformed Lambert $W$-based Colebrook equation in the form more suitable for computing codes used in various engineering software. The best version of the presented explicit approximation gives the value of flow friction factor $f$, for which the Colebrook equation is in balance with the relative error of no more than 0.0096%. Such accuracy achieved without using a large number of computationally expensive logarithmic functions (or non-integer powers) is highly computationally efficient. As reported by Clamond [10], Winning and Coole [11], Biberg [4], Vatankhah [12], etc., functions such as logarithms and non-integer powers require special algorithms with execution of many more floating-point operations compared with the basic arithmetic operations (+,-,*,/) that are executed directly in the Central Processor Unit (CPU) of computers. Apparently, this is the first highly accurate explicit approximation of the Colebrook equation that contains only two computationally expensive functions (two logarithms or as an alternative two functions with non-integer powers) or even less if a combination of Padé approximations [13,14] and symbolic regression is used for a further reduction of the computational burden (where as a result one of the logarithms is approximated by simple rational functions with moderate increase of the maximal relative error).

## 2. Proposed explicit approximations and comparative analysis

The Colebrook equation in the term of the Lambert $W$-function was apparently first proposed in 2018 by Keady [3]. However, as confirmed by Sonnad and Goudar [6] and Brkić [7], the term $W(e^x)$ grows so fast that cannot be evaluated easily even in registers of modern computer due to the overflow runtime error for a certain number of combinations of the input parameters; the Reynolds number $R$ and the relative roughness of inner pipe surface $\varepsilon^*$; where parameter $x$ of equation (2) depends directly on them. The here shown procedure replaces this fast growing term by the much more numerically stable Wright $\Omega$-function [15]. As noted by Lawrence et al. [15], the Wright $\Omega$-function was studied implicitly, without being named, by Wright [16], and named and defined by Corless and Jeffrey [17].

Further about the Colebrook equation transformed in explicit form in term of the Lambert $W$-function can be found in Keady [3], Goudar and Sonnad [18,19], Brkić [20-23], More [24], Sonnad and Goudar [25,26], Clamond [10], Rollmann and Spindler [9], Mikata and Walczak [27], Biberg [4], Vatankhah [12], etc.

## 2.1. Transformation and formulation

The shifted Wright $\Omega$-function transforms the argument $e^x$ to $x$ in the series $W(e^x) = \Omega(x) \approx ln(e^x) - ln(ln(e^x)) + \frac{ln(ln(e^x))}{ln(e^x)}$; where $x = ln(e^x)$. In that way, the undesirably fast growing term $W(e^x)$ in Eq. (2) is approximated accurately through $y \approx x - ln(x) + \frac{ln(x)}{x}$. The transformation is based on unsigned Stirling numbers of the first kind as reported by Rollmann and Spindler [9]. Table 1 shows values of $W(e^x)$ compared with its approximate replacement in the domain of applicability of the Colebrook equation. Without the proposed transformation and simplification, the runtime overflow error occurs during the evaluation of the friction factor $f$ in computers for certain pairs or the Reynolds number $R$ and the relative roughness of inner pipe surface $e^x$; where parameter x of Eq. (2) depends directly on them (#VALUE! is overflow error in Table 1). The values in Table 1 are calculated in MS Excel.

**Table 1.** Values of $W(e^x)$ compared with its approximate replacement $y \approx x - lnx + \frac{lnx}{x}$.

| $W(e^x)$ | $R=4000$ | $R=10^4$ | $R=10^5$ | $R=10^6$ | $R=10^7$ | $R=10^8$ |
|---|---|---|---|---|---|---|
| $\varepsilon^* = 10^{-6}$ | 5.763586714 | 6.552354737 | 8.594740889 | 10.78188015 | 13.94025768 | 26.71930109 |
| $\varepsilon^* = 10^{-5}$ | 5.767379666 | 6.562009418 | 8.694474328 | 11.80401384 | 24.50329461 | 125.7849498 |
| $\varepsilon^* = 10^{-3}$ | 5.805329409 | 6.658658836 | 9.697953496 | 22.29514802 | 124.0554132 | #VALUE! |
| $\varepsilon^* = 10^{-2}$ | 6.186774452 | 7.63459358 | 20.09639172 | 122.325789 | #VALUE! | #VALUE! |
| $\varepsilon^* = 0.05$ | 10.14320931 | 17.90904123 | 120.5960672 | #VALUE! | #VALUE! | #VALUE! |
| $y$ | $R=4000$ | $R=10^4$ | $R=10^5$ | $R=10^6$ | $R=10^7$ | $R=10^8$ |
| $\varepsilon^* = 10^{-6}$ | 5.766606874 | 6.552971455 | 8.592338256 | 10.7784212 | 13.93654591 | 26.71669441 |
| $\varepsilon^* = 10^{-5}$ | 5.770385511 | 6.562602762 | 8.691991603 | 11.80037821 | 24.50049484 | 136.3596559 |
| $\varepsilon^* = 10^{-3}$ | 5.808193728 | 6.659024862 | 9.694862641 | 22.29214094 | 134.073966 | 1246.853296 |
| $\varepsilon^* = 10^{-2}$ | 6.188374207 | 7.633218988 | 20.093168 | 131.7885643 | 1244.552558 | 12371.62215 |
| $\varepsilon^* = 0.05$ | 10.13993873 | 17.90560354 | 129.5034606 | 1242.251823 | 12369.31975 | 123639.9564 |

#VALUE! – Overflow error

The simplifications; $W(e^x) - x \approx ln(x) \cdot \left(\frac{1}{x} - 1\right)$; $\frac{2}{ln10} \approx 0.8686$; $\frac{2 \cdot 2.51}{ln10} \approx 2.18$; and $2.18 \cdot 3.71 \approx 8.0878$, transform the Lambert $W$-based expression of the Colebrook equation in a very accurate explicit approximate form that can be used efficiently in everyday engineering practice; Eq. (3):

$$\left. \begin{array}{c} \frac{1}{\sqrt{f}} \approx 0.8686 \cdot \left[B + ln(B+A) \cdot \left(\frac{1}{B+A} - 1\right)\right] \\ A \approx \frac{R \cdot \varepsilon^*}{8.0878} \\ B \approx ln\left(\frac{R}{2.18}\right) \approx ln(R) - 0.779397488 \end{array} \right\}, \qquad (3)$$

Instead of logarithmic functions in the proposed explicit approximation; Eq. (3), a new form for $B$ and $ln(B + A)$ can be introduced, where $a$ can be any sufficiently large constant, where the larger value of $a$ gives the more accurate approximation of logarithmic function, Eq. (4):

$$\left.\begin{array}{l} B \approx a \cdot \left(\frac{R}{2.18}\right)^{a^{-1}} - a \\ \ln(B + A) \approx a \cdot (B + A)^{a^{-1}} - a \end{array}\right\}, \quad (4)$$

Very accurate results are obtained for $a > 10^5$. Choosing this value, power $a^{-1} = \frac{1}{a}$ is a fraction with integer numerator and denominator, where the appropriate form depends on the programming language and the option with fever floating point operations should be chosen [28].

The forms such as $\left(\frac{R}{2.18}\right)^{0.00001}$ requires evaluation of two transcendental functions because compilers in most programming languages interpret it through $e^{0.00001 \cdot ln\left(\frac{R}{2.18}\right)}$ [10].

For more accurate results $W(e^x) - x \approx \frac{1.038 \cdot \ln(x)}{x+0.332} - \ln(x)$ or $W(e^x) - x \approx \frac{1.0119 \cdot \ln(x)}{x} - \ln(x) + \frac{\ln(x)-2.3849}{x^2}$ can be used. These new approximations were found by symbolic regression software Eureqa [29-31] and they are 2.5 and 16.7 times respectively more accurate compared with the expression $y$ from Table 1. The related approximations are given with Eq. (5) and Eq. (6), respectively.

$$\frac{1}{\sqrt{f}} \approx 0.8686 \cdot \left[B + \frac{1.038 \cdot ln(B+A)}{0.332+B+A} - ln(B+A)\right], \quad (5)$$

$$\frac{1}{\sqrt{f}} \approx 0.8686 \cdot \left[B + \frac{1.0119 \cdot ln(B+A)}{B+A} - ln(B+A) + \frac{ln(B+A)-2.3849}{(B+A)^2}\right], \quad (6)$$

In Eq. (5) and Eq. (6), parameters $A$ and $B$ are the same as in Eq. (3).

## 2.2. Accuracy

With the friction factor $f$ computed using the approximate equations Eq. (3), the Colebrook equation is in balance with the relative error of no more than 0.13%, while using Eq. (5) of no more than 0.045%, and finally, using Eq. (6) of no more than 0.0096%, respectively. Related distribution of errors is shown in Figure 1. The presented approximations require evaluation of only two computationally expensive functions (two logarithms; Eq. (3), Eq. (5) and Eq. (6) or alternatively two non-integer powers; Eq. (4)), and therefore they are not only accurate, but also efficient for calculation.

The here shown approximation; Eq. (3) based on the Wright $\Omega$-function with the relative error up to 0.13% is about ten times more accurate compared with the approximation from Brkić [22], while Eq. (5) more than 25 times and finally, Eq. (6) is more than 100 times more accurate than [22]. The approximations from Brkić [22,23] are based on the Lambert $W$-function.

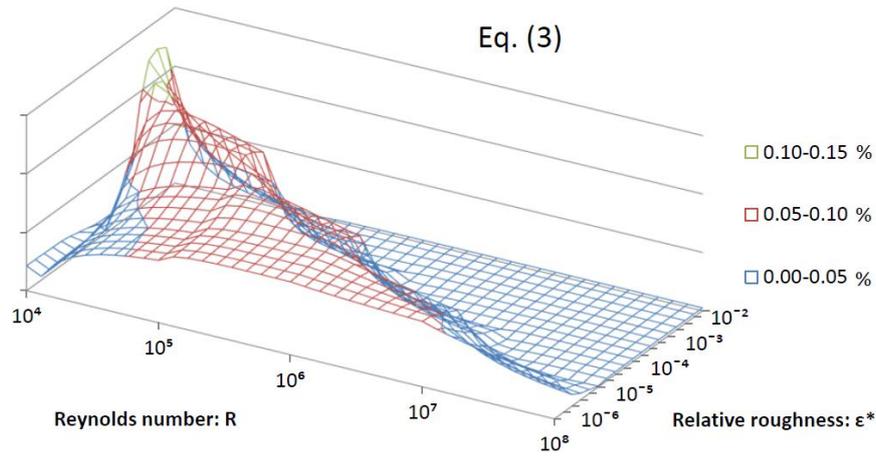

(a)

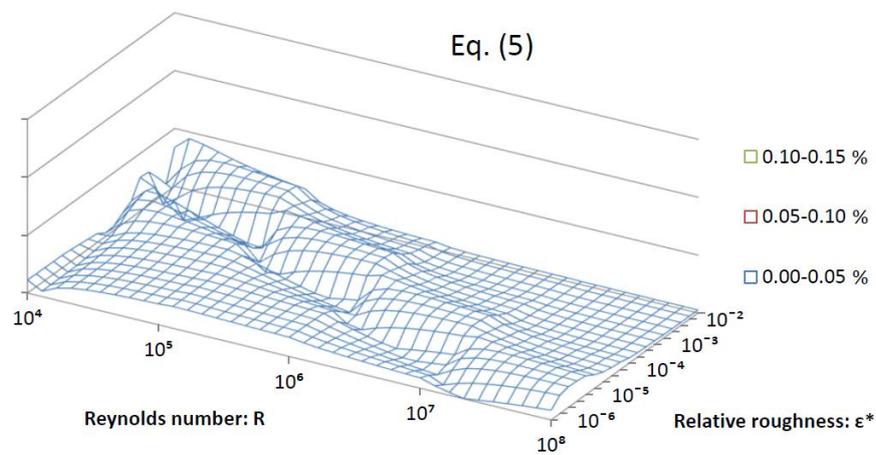

(b)

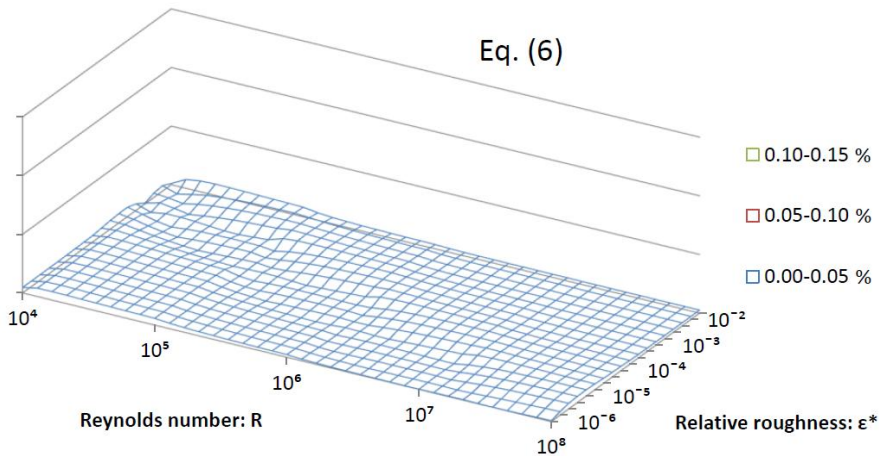

(c)

**Figure 1.** Distribution of the relative error by the proposed explicit approximation of Colebrook's equation; (a): Eq. (3), (b): Eq. (5) and (c): Eq. (6); Comparison

The most accurate approximations available up to date are by Vatankhah [12], Buzzelli [32], Vatankhah and Kouchakzadeh [33], Romeo et al. [34], Zigrang and Sylvester [35] and Serghides [36]. All approximations mentioned or developed in this paper are listed in Appendix while its evaluated relative error in Table 2.

Further about accuracy of explicit approximations to the Colebrook equation can be found in Zigrang and Sylvester [37], Gregory and Fogarasi [38], Brkić [39,40], Winning and Coole [11,41], Brkić and Ćojbašić [42].

*2.3. Complexity and computational burden*

In computer environment, a logarithmic function and non-integer powers require more floating-point operations to be executed in the CPU compared with the simple arithmetic operations such as adding, subtracting, multiplication or division [10-12]. With the relative error of up to 0.0096%, the here proposed explicit approximation of the Colebrook equation; Eq. (6), that contains only two computationally expensive functions, is not only accurate but also sufficiently efficient. Winning and Coole [11] reported relative effort for computation as: Addition-1, Subtraction-1.18, Division-1.35, Multiplication-1.55, Squared-2.18, Square root-2.29, Cubed-2.38, Natural logarithm-2.69, Cubed root-2.71, Fractional exponential-3.32, and Logaritm to base 10-3.37. On the other hand, Biberg [4] adds division in the group of more expensive functions.

For comparison, Table 2 provides number of logarithmic functions and non-integer terms used in available approximations. Table 2 shows only highly accurate approximations with the relative error of no more than 1% according to criterions set by Brkić [29]. All approximations from Table 2 are given in Appendix of this article.

**Table 2.** Number of computationally expensive functions in the available approximations of the Colebrook equations that introduce relative error of no more than 1%.

| Approximation [1] | Maximal relative error % | Function | | Total [2] |
|---|---|---|---|---|
| | | Logarithms | Non-integer powers | |
| Vatankhah | 0.0028% | 1 | 2 | 3(5) |
| Here developed; Eq. (6) | 0.0096% | 2 | 0 | 2 |
| Here developed; Eq. (5) | 0.045%, | 2 | 0 | 2 |
| Here developed; Eq. (3) | 0.13% | 2 | 0 | 2 |
| Here developed; Eq. (4) | 0.13% | 0 | 2 | 2(4) |
| Buzzelli [3] | 0.14% | 2 | 0 | 2 |
| Zigrang and Sylvester | 0.14% | 3 | 0 | 3 |
| Serghides | 0.14% | 3 | 0 | 3 |
| Romeo et al. | 0.14% | 3 | 2 | 5(7) |
| Vatankhah and Kouchakzadeh | 0.15% | 2 | 1 | 3(4) |
| Barr | 0.27% | 2 | 2 | 4(6) |
| Serghides-simple | 0.35% | 2 | 0 | 2 |
| Chen | 0.36% | 2 | 2 | 4(6) |
| Here developed; Eq. (11) | Up to 0.4% | 1 | 0 | 1 |
| Fang et al. | 0.62% | 1 | 3 | 4(7) |
| Papaevangelou et al. | 0.82% | 2 | 1 | 3(4) |

[1] All approximations are listed in Appendix of this paper, [2] in brackets: according to Clamond [10] non-integer powers require evaluation of two computationally expensive functions – logarithm and exponential function, [3] in addition contains also one square root function.

In addition to the here presented, the approximation by Brkić [22,23] is also based on the Lambert W-function, but with four logarithmic functions used, it is also much more computationally expensive (with the relative error of about 2.2% it is also significantly less accurate).

In the next Section, the logarithmic function of $B \approx ln\left(\frac{R}{2.18}\right)$ from Eq. (3) is approximated very accurately trough rational polynomial expression, so complexity and computational additionally decrease. Also, subtraction requires less floating-point operations than division, so the computational cheaper form $B \approx ln(R) - 0.779397488$ should be used instead.

*2.4. Simplifications*

A simple rational approximation of the logarithm term $B$ of the novel Colebrook approximation formulas; Eq. (3), Eq. (5) and Eq. (6), is shown in this Section. The logarithm represents the most computationally expensive operation of the Colebrook formula, while the most its approximations also contain computationally demanded non-integer power terms. In order to reduce computation costs, the idea is to replace the term $B$ that contains the logarithmic function by simple rational functions. A combination of Padé approximation [14,43] and an artificial intelligence symbolic regression procedure [29-31,] is used for this. Although the logarithm is a transcendental function, the found rational approximation remains simple and accurate with the maximal relative error limited to 0.2%. Although this rational approximation of the logarithm is not very nice for a human perception, it is very fast at computers, as it requires only a limited number of basic arithmetic operations to be executed in the CPU.

For the purpose of this simplification, the observed form $B$ from Eq. (3) can be transformed as; Eq. (7):

$$B \approx \ln(R) - \ln(2.18) = \ln(315012.6 \cdot r) - 0.77932 = \ln(r) + \ln(315012.6) - 0.77932 = \ln r + 11.881, \quad (7)$$

In Eq. (7), for term $r = \frac{R}{315012.6}$, constant 315012.6 is carefully selected in order to minimize the error of the /2,3/ order Padé approximation of $\ln(r)$ at the expansion point $r_0 = 1$. The proposed Padé approximant of the /2,3/ order of $\ln(r)$ at point 1 is; Eq. (8):

$$\ln(r) \approx s(r) = \frac{r \cdot (r \cdot (11 \cdot r + 27) - 27) - 11}{r \cdot (r \cdot (3 \cdot r + 27) + 27) + 3}, \quad (8)$$

The value 315012.6 is a weighted average of the Reynolds number $R$ for the turbulent zone valid for the Colebrook equation; $R_{min} = 4000$ and $R_{max} = 10^8$, using the value 0.0063 that was set by numerical experiments in order to minimize the absolute value of the maximum relative error of the Padé approximant of $\ln(r)$ in interval $[R_{min}, R_{max}]$ as $0.0063 \cdot \frac{R_{min}+R_{max}}{2} = 315012.6$. The Padé approximant approximates a certain function very accurately only in a relatively short domain of input parameters. It has been observed that the Padé approximant of $\ln(r)$ at the expansion point $r_0 = 1$ defined by a rational function $s(r)$ approximates $\ln(r)$ with the maximal relative error between -11.8% and 11.8% for all values of the Reynolds number $R$ in the interval $[R_{min}, R_{max}]$. The Padé approximant $s(r)$ has a negligible error for $r \sim 1$, whereas top errors correspond to border points $R_{min}$ and $R_{max}$. For example, for $R_{min} = 4000$, the Padé approximant of $\ln(r)$ is $s\left(\frac{4000}{315012.6}\right) = s(0.012697905) = -3.38744549$ where $R = 315012.6 \cdot r$. Therefore, the value of $\ln(4000)$ is approximated by $\ln(315012.6) + (-3.387445469) = 9.272922448$. The corresponding relative error for $\ln(4000)$, is -11.8%.

Because of $B = \ln\left(\frac{R}{\frac{2 \cdot 2.51}{\ln(10)}}\right) = \ln(R) - \ln\left(\frac{2 \cdot 2.51}{\ln(10)}\right)$ and $\ln(315012.6) - \ln\left(\frac{2 \cdot 2.51}{\ln(10)}\right) \sim 11.881$, the value of $B$ can be approximated as; Eq. (9):

$$B \approx \ln(r) + 11.881 \approx s(r) + 11.881, \tag{9}$$

Further, a symbolic regression technique based on computer software Eureqa [30,31], is used for a further more precise approximation of $\ln(r)$. The aim is to construct a more accurate rational approximation of $\ln(r)$ in comparison with Eq. (9) using two known variables: the ratio $r = \frac{R}{315012.6}$ and its Padé approximation $s(r)$. In order to reduce the burden for the CPU, the symbolic regression model should have a computationally cheap evaluation. For this reason, only rational functions are assumed for the symbolic regression model. To achieve that, 200 carefully selected quasi-random points of $r$ using LPTAU51 algorithm is used [44,45]. For these generated numbers, the Padé approximation $s(r)$ is calculated using Eq. (8). Also, $\ln(r)$ is calculated in order to train the model in Eureqa for the purpose to find a rational approximation of $\ln(r)$ by using $r$ and $s(r)$ pairs. The developed models were successfully tested using 2048 quasi-random points. As a result, value $B$ is approximated by simple rational functions; Eq. (10), with the negligible maximal relative error 0.0765%.

$$B \approx 0.98236 \cdot s + \frac{s^7}{9200.67} + \frac{r}{150.2325} - \frac{r^2}{138187.1651} - \frac{1}{161.124 \cdot r} + 11.881, \tag{10}$$

Here the symbol s denotes the Padé approximant $s(r)$ given by Eq. (8) and $r = \frac{R}{315012.6}$ is its argument.

When the Horner nested representation and the Variable Precision Arithmetic (VPA) at 4 decimal digit accuracy is assumed, the approximation of $B$ can be simplified by Eq. (11):

$$B \approx s \cdot 0.0001086 \cdot s^6 + 0.9824 - \frac{0.006206}{r} - r \cdot 0.000007237 \cdot r - 0.006656 + 11.88, \tag{11}$$

In this case, the maximal relative error remains negligible, 0.0793% compared with $B$ calculated using Eq. (3).

The combined approach with Padé approximant and the symbolic regression introduces in this Section is based on a human observation and introducing the ratio $r = \frac{R}{315012.6}$ with the subsequent symbolic regression of $r$ and $s(r)$ pairs by the Eureqa. The maximal relative error of $B$ introduced by Eq. (11) is small; 0.0793% and in total if it is used instead of $B \approx \ln\left(\frac{R}{2.18}\right)$ from Eq. (3), the total maximal error of the explicit approximation of the Colebrook equation can go up to 0.4%. As can be seen from Table 2, with the only one-log call for $\ln(B + A)$ from Eq. (3), this approximation is the cheapest for computation to date presented extremely accurate explicit approximation of the Colebrook equation.

The here presented combined approach with Padé approximation and the symbolic regression can be also used for faster but still accurate probabilistic modeling of gas networks, which requires millions of model evaluations [46-48].

## 3. Conclusion

Although the implicit Colebrook equation for flow friction is empirical and hence with disputed accuracy, in many cases it is necessary to repeat calculations and to resolve the equation accurately in order to compare scientific results. An iterative solution [49,50] requires extensive computational efforts especially for flow evaluation of complex water or gas pipeline networks [51-53]. Although various available explicit approximations offer a good alternative, they are by the rule very accurate, but too complex and vice versa [39]. In contrary to previous approximations of the Colebrook equation, the here presented relation with the relative error limited to 0.0096% belongs to the group of the most accurate available explicit approximations of the Colebrook equation. Moreover, the here presented approach is also very cheap, as it needs only one or two logarithms (or alternatively two non-integer powers). According to the both criterions; accuracy and complexity, the here presented approximations show interesting performance. For this reason, the here presented approximations can be recommended for implementation in software codes for engineering use.

The Colebrook equation is relevant only for the turbulent flow, while for the full-scale flow different unified equations can be used [54].

**Abbreviations**

The following symbols are used in this paper:

Constants:

$a$ - any $> 10^5$

Variables:

A - variable that depends on $R$ and $\varepsilon^*$ (dimensionless)

B - variable that depends on $R$ (dimensionless)

$f$ – Darcy (Moody) flow friction factor (dimensionless)

$R$ – Reynolds number (dimensionless)

$r$ - variable that depends on $R$ (dimensionless)

$x$ – variable in function on R and $\varepsilon^*$ (dimensionless)

$\varepsilon^*$ – Relative roughness of inner pipe surface (dimensionless)

$\alpha$ – variables defined in Appendix of this paper

Functions:

e – exponential function

ln – natural logarithm

s – Padé approximant

W – Lambert function

Ω – Wright function

**Appendix**

The following explicit approximations of the Colebrook equation are referred in this paper:

-Here developed Eq. (3), Eq. (5), Eq. (6) - Eq. (A.1.1), Eq. (A.1.2), Eq. (A.1.3):

$$\frac{1}{\sqrt{f}} \approx 0.8686 \cdot \left[ B + ln(B+A) \cdot \left( \frac{1}{B+A} - 1 \right) \right], \quad (A.1.1)$$

$$\frac{1}{\sqrt{f}} \approx 0.8686 \cdot \left[ B + \frac{1.038 \cdot ln(B+A)}{0.332+B+A} - ln(B+A) \right], \quad (A.1.2)$$

$$\frac{1}{\sqrt{f}} \approx 0.8686 \cdot \left[ B + \frac{1.0119 \cdot ln(B+A)}{B+A} - ln(B+A) + \frac{ln(B+A)-2.3849}{(B+A)^2} \right], \quad (A.1.3)$$

Where $A \approx \frac{R \cdot \varepsilon^*}{8.0878}$, $B \approx ln\left(\frac{R}{2.18}\right) \approx ln(R) - 0.779397488$

-Here developed Eq. (4); Eq. (A.2.1), Eq. (A.2.2), Eq. (A.2.3):

$$\frac{1}{\sqrt{f}} \approx 0.8686 \cdot \left[ B + \left( a \cdot (B+A)^{a^{-1}} - a \right) \cdot \left( \frac{1}{B+A} - 1 \right) \right], \quad (A.2.1)$$

$$\frac{1}{\sqrt{f}} \approx 0.8686 \cdot \left[ B + \frac{1.038 \cdot \left( a \cdot (B+A)^{a^{-1}} - a \right)}{0.332+B+A} - \left( a \cdot (B+A)^{a^{-1}} - a \right) \right], \quad (A.2.2)$$

$$\frac{1}{\sqrt{f}} \approx 0.8686 \cdot \left[ B + \frac{1.0119 \cdot \left( a \cdot (B+A)^{a^{-1}} - a \right)}{B+A} - \left( a \cdot (B+A)^{a^{-1}} - a \right) + \frac{\left( a \cdot (B+A)^{a^{-1}} - a \right) - 2.3849}{(B+A)^2} \right], \quad (A.2.3)$$

Where $A \approx \frac{R \cdot \varepsilon^*}{8.0878}$, and $B \approx a \cdot \left( \frac{R}{2.18} \right)^{a^{-1}} - a \approx a \cdot (R)^{a^{-1}} - a - 0.779397488$.

As parameter $a$ is larger, the more accurate solution is ($a > 10^5$ gives sufficiently satisfied results).

-Here developed Eq. (11); Eq. (A.3):

Parameter $B$ from the Eqs. (A.1.1)-(A.1.3) and Eqs. (A.2.1)-(A.2.3) should be calculated using Eq. (A.3).

$$\left. \begin{array}{l} B \approx s \cdot (0.0001086 \cdot s^6 + 0.9824) - \frac{0.006206}{r} - r \cdot (0.000007237 \cdot r - 0.006656) + 11.881 \\ r = \frac{R}{315012.6} \\ s \approx s(r) = \frac{r \cdot (r \cdot (11 \cdot r + 27) - 27) - 11}{r \cdot (r \cdot (3 \cdot r + 27) + 27) + 3} \end{array} \right\}, \quad (A.3)$$

-Buzzelli [32]; (A.4):

$$\left. \begin{array}{l} \frac{1}{\sqrt{f}} \approx \alpha_1 - \left( \frac{\alpha_1 + 2 \cdot \log_{10}\left( \frac{\alpha_2}{R} \right)}{1 + \frac{2.18}{\alpha_2}} \right) \\ \alpha_1 \approx \frac{(0.774 \cdot \ln(R)) - 1.41}{1 + 1.32 \cdot \sqrt{\varepsilon^*}} \\ \alpha_2 \approx \frac{\varepsilon^*}{3.7} \cdot R + 2.51 \cdot \alpha_1 \end{array} \right\}, \quad (A.4)$$

-Zigrang and Sylvester [35]; (A.5):

$$\left. \begin{array}{l} \frac{1}{\sqrt{f}} \approx -2 \cdot \log_{10}\left( \frac{\varepsilon^*}{3.7} - \frac{5.02}{R} \cdot \alpha_3 \right) \\ \alpha_3 \approx \log_{10}\left( \frac{\varepsilon^*}{3.7} - \frac{5.02}{R} \cdot \alpha_4 \right) \\ \alpha_4 \approx \log_{10}\left( \frac{\varepsilon^*}{3.7} - \frac{13}{R} \right) \end{array} \right\} \quad (A.5)$$

-Serghides [36]; (A.6):

$$\left.\begin{aligned}
\frac{1}{\sqrt{f}} &\approx \alpha_5 - \frac{(\alpha_6 - \alpha_5)^2}{\alpha_7 - 2 \cdot \alpha_6 + \alpha_5} \\
\alpha_5 &\approx -2 \cdot \log_{10}\left(\frac{\varepsilon^*}{3.7} - \frac{12}{R}\right) \\
\alpha_6 &\approx -2 \cdot \log_{10}\left(\frac{\varepsilon^*}{3.7} - \frac{2.51}{R} \cdot \alpha_5\right) \\
\alpha_7 &\approx -2 \cdot \log_{10}\left(\frac{\varepsilon^*}{3.7} - \frac{2.51}{R} \cdot \alpha_6\right)
\end{aligned}\right\} \quad (A.6)$$

-Romeo et al. [34]; (A.7):

$$\left.\begin{aligned}
\frac{1}{\sqrt{f}} &\approx -2 \cdot \log_{10}\left(\frac{\varepsilon^*}{3.7065} - \frac{5.0272}{R} \cdot \alpha_8\right) \\
\alpha_8 &\approx \log_{10}\left(\frac{\varepsilon^*}{3.827} - \frac{4.567}{R} \cdot \alpha_9\right) \\
\alpha_9 &\approx \log_{10}\left(\left(\frac{\varepsilon^*}{7.7918}\right)^{0.9924} + \left(\frac{5.3326}{208.815 + R}\right)^{0.9345}\right)
\end{aligned}\right\} \quad (A.7)$$

-Vatankhah and Kouchakzadeh [33]; (A.8):

$$\left.\begin{aligned}
\frac{1}{\sqrt{f}} &\approx 0.8686 \cdot \ln\left(\frac{0.4587 \cdot R}{(\alpha_{10} - 0.31)^{\alpha_{11}}}\right) \\
\alpha_{10} &\approx 0.124 \cdot R \cdot \varepsilon^* + \ln(0.1587 \cdot R) \\
\alpha_{11} &\approx \frac{\alpha_{10}}{\alpha_{10} + 0.9633}
\end{aligned}\right\} \quad (A.8)$$

-Barr [55]; (A.9):

$$\left.\begin{aligned}
\frac{1}{\sqrt{f}} &\approx -2 \cdot \log_{10}\left(\frac{\varepsilon^*}{3.7} + \frac{4.518 \cdot \log_{10}\left(\frac{R}{7}\right)}{\alpha_{12}}\right) \\
\alpha_{12} &\approx R \cdot \left(1 + \frac{R^{0.52}}{29} \cdot (\varepsilon^*)^{0.7}\right)
\end{aligned}\right\} \quad (A.9)$$

-Serghides-simple [36]; (A.10):

$$\left.\begin{aligned}
\frac{1}{\sqrt{f}} &\approx 4.781 - \frac{(\alpha_{13} - 4.781)^2}{\alpha_{14} - 2 \cdot \alpha_{13} + 4.781} \\
\alpha_{13} &\approx -2 \cdot \log_{10}\left(\frac{\varepsilon^*}{3.7} - \frac{12}{R}\right) \\
\alpha_{14} &\approx -2 \cdot \log_{10}\left(\frac{\varepsilon^*}{3.7} - \frac{2.51}{R} \cdot \alpha_{13}\right)
\end{aligned}\right\} \quad (A.10)$$

-Chen [56]; (A.11):

$$\left.\begin{aligned}
\frac{1}{\sqrt{f}} &\approx -2 \cdot \log_{10}\left(\frac{\varepsilon^*}{3.7065} - \frac{5.0452}{R} \cdot \alpha_{15}\right) \\
\alpha_{15} &\approx \log_{10}\left(\frac{(\varepsilon^*)^{1.1098}}{2.8257} + \frac{5.8506}{R^{0.8981}}\right)
\end{aligned}\right\} \quad (A.11)$$

-Fang et al. [57]; (A.12):

$$\left.\begin{aligned}
\frac{1}{\sqrt{f}} &\approx \left(1.613 \cdot \left(\ln\left(0.234 \cdot (\varepsilon^*)^{1.1007} - \alpha_{16}\right)\right)^{-2}\right)^{-2} \\
\alpha_{16} &\approx \frac{60.525}{R^{1.1105}} + \frac{56.291}{R^{1.0712}}
\end{aligned}\right\} \quad (A.12)$$

-Papaevangelou et al. [58]; (A.13):

$$\frac{1}{\sqrt{f}} \approx \left(\frac{0.2479 - 0.0000947 \cdot (7 - \log_{10}(R))^4}{\left(\log_{10}\left(\frac{\varepsilon^*}{3.615} + \frac{7.366}{R^{0.9142}}\right)\right)^2}\right)^{-2} \quad (A.13)$$

-Vatankhah.[12]; (A.14):

$$\left. \begin{array}{l} \frac{1}{\sqrt{f}} \approx 0.8686 \cdot \ln\left(\frac{0.3984 \cdot R}{(0.8686 \cdot \alpha_{17})^{\frac{\alpha_{17}}{\alpha_{17}+\alpha_{18}}}}\right) \\ \alpha_{17} \approx 0.12363 \cdot R \cdot \varepsilon^* + \ln(0.3984 \cdot R) \\ \alpha_{18} \approx 1 + \frac{1}{\frac{1+\alpha_{17}}{0.5 \cdot \ln(0.8686 \cdot \alpha_{17})} - \frac{1+4 \cdot \alpha_{17}}{3 \cdot (1+\alpha_{17})}} \end{array} \right\} \quad (A.14)$$